\magnification 1200


\input plainenc
\input amssym
\fontencoding{T2A}
\inputencoding{utf-8}
\tolerance 7500
\binoppenalty 10000
\relpenalty 10000
\parindent 1.5em

\catcode`\@ 11
\font\TITLE labx1440
\font\tenrm larm1000
\font\cmtenrm cmr10
\font\tenit lati1000
\font\tenbf labx1000
\font\tentt latt1000
\font\teni cmmi10 \skewchar\teni '177
\font\tensy cmsy10 \skewchar\tensy '60
\font\tenex cmex10
\font\teneufm eufm10
\font\tenmsb msbm10
\font\ninerm larm0900

\font\ninei cmmi9 \skewchar\ninei '177
\font\ninesy cmsy9 \skewchar\ninesy '60

\font\eightrm larm0800
\font\cmeightrm cmr8
\font\eightit lati0800
\font\eightbf labx0800
\font\eighti cmmi8 \skewchar\eighti '177
\font\eightsy cmsy8 \skewchar\eightsy '60
\font\eightex cmex8
\font\eighteufm eufm8
\font\eightmsb msbm8

\font\cmsevenrm cmr7

\font\sevenbf labx0700

\font\seveni cmmi7 \skewchar\seveni '177
\font\sevensy cmsy7 \skewchar\sevensy '60

\font\seveneufm eufm7
\font\sevenmsb msbm7

\font\sixi cmmi6 \skewchar\sixi '177
\font\sixsy cmsy6 \skewchar\sixsy '60

\font\cmfiverm cmr5
\font\fivebf labx0500
\font\fivei cmmi5 \skewchar\fivei '177
\font\fivesy cmsy5 \skewchar\fivesy '60
\font\fiveeufm eufm5
\font\fivemsb msbm5
\font\tencmmib cmmib10 \skewchar\tencmmib '177
\font\ninecmmib cmmib9 \skewchar\ninecmmib '177
\font\eightcmmib cmmib8 \skewchar\eightcmmib '177
\font\sevencmmib cmmib7 \skewchar\sevencmmib '177
\font\sixcmmib cmmib6 \skewchar\sixcmmib '177
\font\fivecmmib cmmib5 \skewchar\fivecmmib '177
\newfam\cmmibfam
\textfont\cmmibfam\tencmmib \scriptfont\cmmibfam\sevencmmib
\scriptscriptfont\cmmibfam\fivecmmib
\def\tenpoint{\def\rm{\fam0\tenrm}\def\it{\fam\itfam\tenit}%
	\def\bf{\fam\bffam\tenbf}\def\tt{\fam\ttfam\tentt}%
	\textfont0\cmtenrm \scriptfont0\cmsevenrm \scriptscriptfont0\cmfiverm%
	\textfont1\teni \scriptfont1\seveni \scriptscriptfont1\fivei%
	\textfont2\tensy \scriptfont2\sevensy \scriptscriptfont2\fivesy%
	\textfont3\tenex \scriptfont3\tenex \scriptscriptfont3\tenex%
	\textfont\itfam\tenit%
	\textfont\bffam\tenbf \scriptfont\bffam\sevenbf%
	\scriptscriptfont\bffam\fivebf%
	\textfont\eufmfam\teneufm \scriptfont\eufmfam\seveneufm%
	\scriptscriptfont\eufmfam\fiveeufm%
	\textfont\msbfam\tenmsb\scriptfont\msbfam\sevenmsb%
	\scriptscriptfont\msbfam\fivemsb%
	\textfont\cmmibfam\tencmmib \scriptfont\cmmibfam\sevencmmib%
	\scriptscriptfont\cmmibfam\fivecmmib%
	\normalbaselineskip 12pt%
	\setbox\strutbox\hbox{\vrule height8.5pt depth3.5pt width0pt}%
	\normalbaselines\rm}

\def\eightpoint{\def\rm{\fam 0\eightrm}\def\it{\fam\itfam\eightit}%
	\def\bf{\fam\bffam\eightbf}%
	\textfont0\cmeightrm \scriptfont0\cmfiverm \scriptscriptfont0\cmfiverm%
	\textfont1\eighti \scriptfont1\fivei \scriptscriptfont1\fivei%
	\textfont2\eightsy \scriptfont2\fivesy \scriptscriptfont2\fivesy%
	\textfont3\eightex \scriptfont3\eightex \scriptscriptfont3\eightex%
	\textfont\itfam\eightit%
	\textfont\bffam\eightbf \scriptfont\bffam\fivebf%
	\scriptscriptfont\bffam\fivebf%
	\textfont\eufmfam\eighteufm \scriptfont\eufmfam\fiveeufm%
	\scriptscriptfont\eufmfam\fiveeufm%
	\textfont\msbfam\eightmsb\scriptfont\msbfam\fivemsb%
	\scriptscriptfont\msbfam\fivemsb%
	\textfont\cmmibfam\eightcmmib \scriptfont\cmmibfam\fivecmmib%
	\scriptscriptfont\cmmibfam\fivecmmib%
	\normalbaselineskip 11pt%
	\abovedisplayskip 5pt%
	\belowdisplayskip 5pt%
	\setbox\strutbox\hbox{\vrule height7pt depth2pt width0pt}%
	\normalbaselines\rm}

\def\No{\char 157}
\def\empty{}

\catcode`\" 12
\def\k@vy{"}
\def\bolddelta{\mathchar"0\hexnumber@\cmmibfam 0E}

\def\dot{\mathaccent"705F }

\catcode`\" 13
\def"#1{\ifx#1<\char 190\relax\else\ifx#1>\char 191\relax\else%
	\ifx#1`\char 16\relax\else\ifx#1'\char 17\relax\else #1\fi\fi\fi\fi}

\long\def\hyp@rprov#1#2#3#4{\hbox{\ifx\hyperrefs\undefined #1\else%
	\ifx#3\relax\else\fi #1%
	\ifx#3\relax\else\fi\fi}}

\newcount\c@section
\newcount\c@subsection
\newcount\c@subsubsection
\newcount\c@equation
\newcount\c@bibl
\newcount\c@enum
\c@section 0
\c@subsection 0
\c@subsubsection 0
\c@equation 0
\c@bibl 0
\newdimen\d@enum
\d@enum 0pt
\def\lab@l{}
\def\label#1{\immediate\write 11{\string\newl@bel{#1}{\lab@l}}\ifhmode\unskip\fi}

\def\bibitem#1{\global\advance\c@bibl 1\par\noindent{\hyp@rprov{%
	\ninerm [\number\c@bibl]}{name=\k@vy cite:\number\c@bibl\k@vy}\relax{}~}%
	\immediate\write 11{\string\newbl@bel{#1}{\number\c@bibl}}}

\def\section#1{\global\advance\c@section 1
	{\par\vskip 3ex plus 0.5ex minus 0.1ex
	\rightskip 0pt plus 1fill\leftskip 0pt plus 1fill\noindent
	\hyp@rprov{{\bf\S\thinspace\number\c@section .~#1}}{%
	name=\k@vy sect:\number\c@section\k@vy}%
	\relax{}\par\penalty 10000\vskip 1ex plus 0.25ex}
	\gdef\lab@l{\number\c@section.:\number\pageno}
	\c@subsection 0
	\c@subsubsection 0
	\c@equation 0
}
\def\subsection{\global\advance\c@subsection 1
	\par\vskip 1ex plus 0.1ex minus 0.05ex\indent
	\hyp@rprov{{\bf\number\c@subsection.{}}}{%
	name=\k@vy sect:\number\c@section.\number\c@subsection\k@vy}\relax{}%
	\gdef\lab@l{\number\c@section.\number\c@subsection:\number\pageno}
	\c@subsubsection 0
}
\def\subsubsection{\global\advance\c@subsubsection 1
	\par\vskip 1ex plus 0.1ex minus 0.05ex\indent%
	\hyp@rprov{{\bf\number\c@subsection.\number\c@subsubsection.{}}}{%
	name=\k@vy sect:\number\c@section.\number\c@subsection.%
	\number\c@subsubsection\k@vy}\relax{}%
	\gdef\lab@l{\number\c@section.\number\c@subsection.%
		\number\c@subsubsection:\number\pageno}
}
\def\equation{\global\advance\c@equation 1
	\gdef\lab@l{\number\c@section.\number\c@equation:\number\pageno}
	\hyp@rprov{{\rm\number\c@equation}}{name=\k@vy eq:\number\c@section(%
	\number\c@equation)\k@vy}\relax{}%
}
\catcode`\# 11
\def\sh@rp{#}
\catcode`\# 6
\def\ref@ref#1.#2:#3:{\def\REF@{#2}\ifx\REF@\empty%
	\hyp@rprov{{\rm \S\thinspace#1}}{href=\k@vy\sh@rp sect:#1\k@vy}0{rgb 0 0 1}%
	\else\hyp@rprov{\ifnum #1=\c@section {\rm #2}%
	\else {\rm \S\thinspace#1.#2}\fi}{href=\k@vy\sh@rp sect:#1.#2\k@vy}0{rgb 0 0 1}\fi
}
\def\ref@pageref#1:#2:{#2}
\def\ref@eqref#1.#2:#3:{\hyp@rprov{\ifnum #1=\c@section {\rm (#2)}\else%
	{\rm \S\thinspace#1$\,$(#2)}\fi}{href=\k@vy\sh@rp eq:#1(#2)\k@vy}0{rgb 0 0 1}%
}
\def\ref#1{\expandafter\ifx\csname l@#1\endcsname\relax
	{\bf ??}\message{^^J Reference #1 undefined!^^J}%
	\else\edef\mur@{\csname l@#1\endcsname :}%
	{\expandafter\ref@ref\mur@}\fi}
\def\pageref#1{\expandafter\ifx\csname l@#1\endcsname\relax
	{\bf ??}\message{^^J Reference #1 undefined!^^J}%
	\else\edef\mur@{\csname l@#1\endcsname :}%
	{\expandafter\ref@pageref\mur@}\fi}
\def\eqref#1{\expandafter\ifx\csname l@#1\endcsname\relax
	{(\bf ??)}\message{^^J Reference (#1) undefined!^^J}%
	\else\edef\mur@{\csname l@#1\endcsname :}%
	{\expandafter\ref@eqref\mur@}\fi}
\def\cite#1{\expandafter\ifx\csname bl@#1\endcsname\relax
	{\bf ??}\message{^^J Citation #1 undefined!^^J}%
	\else\hyp@rprov{{\bf\csname bl@#1\endcsname}}{href=\k@vy\sh@rp cite:%
	\expandafter\number\csname bl@#1\endcsname\k@vy}0{rgb 0 0 1}\fi}

\long\def\enumerate#1#2{%
	\setbox0\hbox{$#1^{\circ}.\ $}\d@enum\wd0\global\advance\d@enum 2pt\c@enum 0%
	{\def\item{\global\advance\c@enum 1\par\hskip 0pt%
	\hbox to \d@enum{$\number\c@enum^{\circ}$.\hss}}%
	\par\smallskip #2\par\smallskip}
}
\def\Wo{{\mathpalette\Wo@{}}W}
\def\Wo@#1{\setbox0\hbox{$#1 W$}\dimen@\ht0\dimen@ii\wd0\raise0.65\dimen@%
\rlap{\kern0.35\dimen@ii$#1{}^\circ$}}
\catcode`\@ 12

\def\proof{\par\medskip{\rm Д$\,$о$\,$к$\,$а$\,$з$\,$а$\,$т$\,$е$\,$%
	л$\,$ь$\,$с$\,$т$\,$в$\,$о.}\ }
\def\endproof{{\parfillskip 0pt\hfill$\square$\par}\medskip}

\def\newl@bel#1#2{\expandafter\gdef\csname l@#1\endcsname{#2}}
\def\newbl@bel#1#2{\expandafter\gdef\csname bl@#1\endcsname{#2}}
\def\contentsline#1#2#3{\expandafter\gdef\csname contents@\number\c@content%
	\endcsname{\par\ifnum #1=0\vskip 4pt\else\vskip 1pt\fi\noindent%
	\hbox to \hsize{\ifnum #1=0{\bf #2\thinspace}\else\hskip 0.5cm%
	#2\thinspace\fi\leaders\hbox to 0.25cm{.}\hfill #3}}%
	\global\advance\c@content 1
}
\openin 11=\jobname .aux
\ifeof 11
	\closein 11\relax
\else
	\closein 11
	\input\jobname .aux
	\relax
\fi

\immediate\openout 11=\jobname.aux


\hsize 17.5truecm
\vsize 25truecm
\hoffset -0.25truecm
\voffset -0.75truecm
\tolerance 9999

\def\Re{\mathop{\rm Re}}

\def\im{\mathop{\rm im}}
\def\Sp{\mathop{\rm Sp}}
\def\Arg{\mathop{\rm Arg}}

\parindent 1em
\tenpoint\frenchspacing
\leftline{\ }\vskip 0.25cm
{\leftskip 0cm plus 1fill\rightskip 0cm plus 1fill\parindent 0cm\baselineskip 15pt
\TITLE Метод диагонализации в теории асимптотических оценок решений систем
обыкновенных дифференциальных уравнений
\par\vskip 0.15cm\rm А.$\,$А.~Владимиров\par}
\vskip 0.25cm
$$
	\vbox{\hsize 0.75\hsize\leftskip 0cm\rightskip 0cm
	\eightpoint\rm
	{\bf Аннотация:\/} В статье рассматривается вопрос об асимптотическом
	поведении решений систем линейных обыкновенных дифференциальных уравнений
	$\dot x=A_\nu x+f$ при больших значениях параметра $\nu\in\frak A$.
	Излагается подход, позволяющий сводить широкий класс задач к тривиально
	решаемым задачам с диагональной матрицей коэффициентов. Этот подход основан
	на приближении суммируемых матричнозначных функций $A_\nu$ некоторыми
	специально сконструированными гладкими.\par
	}
$$

\vskip 0.25cm
\section{Введение}
\subsection
В последнее время возобновился \hbox{[\cite{SSh:2020}, \cite{Sh:2021}]} давний
интерес \hbox{[\cite{Na:1969}: \S$\,$4]} к проблеме изучения асимптотического
поведения решений граничных задач для уравнений вида $\dot x=A_\nu x+f$ при больших
значениях параметра $\nu\in\frak A$. Здесь $\frak A$ есть произвольно фиксированное
направленное множество; в типичных для конкретных задач ситуациях в роли таких
множеств выступают секторы комплексной плоскости. Целью настоящей статьи является
изложение конструкции, дающей достаточно простой общий подход к изучению такого
рода задач.

\subsection
Далее мы будем использовать следующую систему понятий и обозначений.
Символом $\frak H$ мы обозначим некоторое конечномерное унитарное пространство
положительной размерности. Символом $\frak B$ мы обозначим пространство действующих
в $\frak H$ операторов, снабжённое стандартной спектральной нормой
$$
	\|A\|_{\frak B}\rightleftharpoons\sup_{x\in\frak H\;:\;\|x\|_\frak H=1}
		\|Ax\|_{\frak H}.
$$
Наконец, мы будем рассматривать функциональные пространства $C([0,1];\,\frak H)$,
$C([0,1];\,\frak B)$, $L_1([0,1];\,\frak H)$, $L_1([0,1];\,\frak B)$,
$W_1^1([0,1];\,\frak H)$ и $W_1^1([0,1];\,\frak B)$ с нормами
$$
	\displaylines{\|x\|_{C([0,1];\,\frak H)}\rightleftharpoons
		\sup_{t\in [0,1]}\|x(t)\|_{\frak H},\quad
		\|A\|_{C([0,1];\,\frak B)}\rightleftharpoons
		\sup_{t\in [0,1]}\|A(t)\|_{\frak B},\cr
		\|x\|_{L_1([0,1];\,\frak H)}\rightleftharpoons
		\int_0^1\|x(t)\|_{\frak H}\,dt,\quad
		\|A\|_{L_1([0,1];\,\frak B)}\rightleftharpoons
		\int_0^1\|A(t)\|_{\frak B}\,dt,\cr
		\|x\|_{W_1^1([0,1];\,\frak H)}\rightleftharpoons
		\int_0^1\bigl\{\|x(t)\|_{\frak H}+
			\|\dot x(t)\|_{\frak H}\bigr\}\,dt,\quad
		\|A\|_{W_1^1([0,1];\,\frak B)}\rightleftharpoons
		\int_0^1\bigl\{\|A(t)\|_{\frak B}+
			\|\dot A(t)\|_{\frak B}\bigr\}\,dt.
	}
$$
Символом $\Wo_1^1([0,1];\,\frak B)$ мы, как обычно, будем обозначать
подпространство внутри $W_1^1([0,1];\,\frak B)$, выделенное условием $A(0)=A(1)=0$.
Символом $\Sp A$ мы, как обычно, будем обозначать след оператора $A\in\frak B$.
В частности, след $\Sp 1$ тождественного оператора совпадает с размерностью
пространства $\frak H$. Кроме того, символом $\omega(A)$ мы будем обозначать правую
границу числовой области значений оператора $A\in\frak B$, то есть величину вида
$$
	\omega(A)\rightleftharpoons\sup_{x\in\frak H\;:\;\|x\|_\frak H=1}
		\Re\langle Ax,x\rangle_{\frak H}.
$$
Легко видеть, что для всякой оператор-функции $A\in L_1([0,1];\,\frak B)$
почти всюду определённые композиции $\Sp\circ A$ и $\omega\circ A$ являются
суммируемыми числовыми функциями.

\subsection
Простейшая ситуация возникает, когда коэффициенты изучаемого уравнения представляют
собой "<малое"> возмущение коэффициентов некоторого базового уравнения.
Здесь справедливо следующее хорошо известное несложное утверждение.

\subsubsection\label{prop:0}
{\it Пусть $P\in\frak B$~--- некоторый ортопроектор, а оператор-функции
$A,V\in L_1([0,1];\,\frak B)$ вместе с парой величин $\theta,\gamma\in (0,1)$
удовлетворяют тождеству $A(t)P\equiv PA(t)$, оценке $\|V\|_{L_1([0,1];\,\frak B)}
\leqslant\theta\gamma$ и не зависящей от выбора отрезка $[a,b]\subseteq [0,1]$
оценке
$$
	\int_a^b\omega((2P-1)A(t))\,dt\leqslant-\ln\gamma.\leqno(\equation)
$$\label{eq:00}%
Тогда решения граничных задач
$$
	\displaylines{\dot x=(A+V)x+f,\quad Px(0)+(1-P)x(1)=\xi,\cr
		\dot x_0=Ax_0+f,\quad Px_0(0)+(1-P)x_0(1)=\xi,
	}
$$
где $f\in L_1([0,1];\,\frak H)$ и $\xi\in\frak H$ зафиксированы произвольным
образом, однозначно определены и подчиняются оценке
$$
	\|x-x_0\|_{C([0,1];\,\frak H)}\leqslant{\theta\over 1-\theta}\,
		\|x_0\|_{C([0,1];\,\frak H)}.
$$
}%
\proof
Существование и однозначная определённость вектор-функции
$x_0\in W_1^1([0,1];\,\frak H)$ немедленно вытекает из известных (см., например,
\hbox{[\cite{Na:1969}: \S$\,$16, Теорема~1]}) фактов о разрешимости начальных
задач для линейных уравнений с суммируемыми коэффициентами. Искомая вектор-функция
$x\in W_1^1([0,1];\,\frak H)$ при этом представляет собой в точности решение
интегрального уравнения
$$
	x(t)=x_0(t)+\int_0^1 [P-\vartheta(s-t)]\cdot M(t)[M(s)]^{-1}V(s)x(s)\,ds,
	\leqno(\equation)
$$\label{eq:0}%
где $\vartheta\in L_\infty(\Bbb R)$ есть функция Хэвисайда, а оператор-функция
$M\in W_1^1([0,1];\,\frak B)$ представляет собой решение начальной задачи
$\dot M=AM$, $M(0)=1$. Предположение~\eqref{eq:00} влечёт при этом справедливость
оценок
$$
	\bigl\|[P-\vartheta(s-t)]\cdot M(t)[M(s)]^{-1}\bigr\|_{\frak B}\leqslant
		\gamma^{-1},
$$
означающих, что~\eqref{eq:0} есть уравнение неподвижной точки некоторого
действующего в пространстве $C([0,1];\,\frak H)$ сжимающего отображения
с коэффициентом сжатия $\theta$.
\endproof

Сфера непосредственного действия утверждения~\ref{prop:0} является достаточно
узкой. Однако хорошо известно \hbox{[\cite{P:1974}: \S$\,$17.Л]},
что вектор-функция $y\in W_1^1([0,1];\,\frak H)$, получаемая из решения
уравнения $\dot x=Ax+f$ посредством подстановки $x=Sy$, где оператор-функция
$S\in W_1^1([0,1];\,\frak B)$ при всех $t\in [0,1]$ принимает обратимые значения,
подчиняется уравнению
$$
	\dot y=[S^{-1}AS-S^{-1}\dot S]y+S^{-1}f\leqno(\equation)
$$\label{eq:300}%
с достаточно сильно изменённой матрицей коэффициентов. Во многих типичных случаях
оказывается, что путём замены, конструируемой некоторым стандартным образом,
преобразованное уравнение~\eqref{eq:300} окажется заведомо удовлетворяющим условиям
утверждения~\ref{prop:0} относительно некоторой простой базовой системы. Описанию
именно такой конструкции преобразования уравнения и посвящается дальнейшее.


\section{Гладкие каркасы}
\subsection
Пусть $\Pi\rightleftharpoons\{P_0,P_-,P_+\}$ есть тройка удовлетворяющих равенству
$P_0+P_-+P_+=1$ ортопроекторов на попарно ортогональные подпространства
пространства $\frak H$. В дальнейшем
{\it\hbox{$\Pi$}-ди\-а\-го\-на\-ли\-за\-ци\-ей\/} оператора $A\in\frak B$ мы будем
называть оператор вида
$$
	\Delta_\Pi(A)\rightleftharpoons P_0AP_0+P_-AP_-+P_+AP_+,
$$
заведомо подчинённый оценке $\|\Delta_\Pi(A)\|_{\frak B}\leqslant\|A\|_{\frak B}$.
{\it\hbox{$\Pi$-кар}\-касом\/} мы также будем называть пару $\{B,C\}$
оператор-функций $B\in W_1^1([0,1];\,\frak B)$ и $C\in\Wo_1^1([0,1];\,\frak B)$,
удовлетворяющую следующим двум условиям:
\enumerate{9}{%
\item При любом выборе значения $t\in [0,1]$ оператор $B(t)$ нормален
и перестановочен с каждым из проекторов $P_0$, $P_-$ и $P_+$.
\item Пусть $d_\Pi(B)$ есть наибольшая из величин, равномерно по $t\in [0,1]$
оценивающих снизу взаимные расстояния между частями $\Sigma_0(t)$, $\Sigma_-(t)$
и $\Sigma_+(t)$ спектра оператора $B(t)$, отвечающими разложению пространства
$\frak H$ в прямую сумму инвариантных подпространств
$$
	\frak H=\im P_0\oplus\im P_-\oplus\im P_+.
$$
Тогда справедлива оценка
$$
	\|C\|_{C([0,1];\,\frak B)}<{d_\Pi(B)\over 8\cdot\Sp 1}.\leqno(\equation)
$$\label{eq:99}%
}

\medskip
Формулировка определения \hbox{$\Pi$-кар}\-ка\-са дана нами в предположении,
что ни один из проекторов $P_-$ и $P_+$ не обращается в нуль. Видоизменение этого
определения (и дальнейших связанных с ним рассуждений) для случая, когда такое
обращение имеет место, носит тривиальный характер, и мы не станем на нём специально
останавливаться.

\subsection
Основное свойство \hbox{$\Pi$-кар}\-ка\-сов даётся следующим утверждением.

\subsubsection\label{prop:kar}
{\it Для всякого \hbox{$\Pi$-кар}\-ка\-са $\{B,C\}$ существует оператор-функция
$S\in W_1^1([0,1];\,\frak B)$, удовлетворяющая равенствам $S(0)=S(1)=1$ и оценкам
$$
	\|S(t)-1\|_{\frak B}\leqslant{4\cdot\Sp 1\over d_\Pi(B)}\cdot
		\|C(t)\|_{\frak B},\leqno(\equation)
$$\label{eq:111}%
$$
	\abovedisplayskip -1\jot
	\|[S(t)]^{-1}-1\|_{\frak B}\leqslant{8\cdot\Sp 1\over d_\Pi(B)}\cdot
		\|C(t)\|_{\frak B},\leqno(\equation)
$$\label{eq:112}%
$$
	\abovedisplayskip -1\jot
	\|\dot S\|_{L_1([0,1];\,\frak B)}\leqslant {4\cdot\Sp 1\over
		d_\Pi^2(B)}\cdot\varkappa(B,C),\leqno(\equation)
$$\label{eq:113}%
$$
	\abovedisplayskip -1\jot
	\biggl\|[S(t)]^{-1}\,[B(t)+C(t)]\,S(t)-\Delta_\Pi[B(t)+C(t)]
		\biggr\|_{\frak B}\leqslant
		{8\cdot\Sp 1\over d_\Pi(B)}\cdot\|C(t)\|_{\frak B}^2,
	\leqno(\equation)
$$\label{eq:114}%
где положено
$$
	\varkappa(B,C)\rightleftharpoons\int_0^1\biggl[6\,\|C(t)\|_{\frak B}
		\cdot\|\dot B(t)\|_{\frak B}+\bigl[4\,\|C(t)\|_{\frak B}+
		d_\Pi(B)\bigr]\cdot\|\dot C(t)\|_{\frak B}\biggr]\,dt.
$$
}%
\proof
Обозначим символами $\Gamma_0(t)$, $\Gamma_-(t)$ и $\Gamma_+(t)$ положительно
ориентированные контуры, ограничивающие \hbox{$(d_\Pi(B)/2)$-ок}\-рест\-но\-сти
указанных ранее частей $\Sigma_0(t)$, $\Sigma_-(t)$ и $\Sigma_+(t)$ спектра
оператора $B(t)$. Длины этих контуров не превосходят величин $\pi d_\Pi(B)\cdot
\Sp P_0$, $\pi d_\Pi(B)\cdot\Sp P_-$ и $\pi d_\Pi(B)\cdot\Sp P_+$, соответственно.
При этом очевидным образом имеют место равенства
$$
	\leqalignno{P_0&=-{1\over 2\pi i}\int_{\Gamma_0(t)}(B(t)-z)^{-1}\,dz,\cr
	P_-&=-{1\over 2\pi i}\int_{\Gamma_-(t)}(B(t)-z)^{-1}\,dz,\cr
	P_+&=-{1\over 2\pi i}\int_{\Gamma_+(t)}(B(t)-z)^{-1}\,dz.
	}
$$
Ввиду выполнения оценки~\eqref{eq:99} могут быть построены также следующие
семейства проекторов:
$$
	\leqalignno{Q_0(t)&=-{1\over 2\pi i}\int_{\Gamma_0(t)}
		(B(t)+C(t)-z)^{-1}\,dz,\cr
	Q_-(t)&=-{1\over 2\pi i}\int_{\Gamma_-(t)}(B(t)+C(t)-z)^{-1}\,dz,\cr
	Q_+(t)&=-{1\over 2\pi i}\int_{\Gamma_+(t)}(B(t)+C(t)-z)^{-1}\,dz.
	}
$$
При этом для оператор-функции
$$
	t\mapsto Q_0(t)-P_0={1\over 2\pi i}\int_{\Gamma_0(t)}
		(B(t)+C(t)-z)^{-1}\cdot C(t)\cdot (B(t)-z)^{-1}\,dz
$$
справедливы оценки
$$
	\|Q_0(t)-P_0\|_{\frak B}\leqslant {4\cdot\Sp P_0\over d_\Pi(B)}\cdot
		\|C(t)\|_{\frak B},
$$
а также оценки
$$
	\|\dot Q_0\|_{L_1([0,1];\,\frak B)}\leqslant {4\cdot\Sp P_0\over
		d_\Pi^2(B)}\cdot\varkappa(B,C).
$$
Аналогичным образом устанавливается справедливость оценок
$$
	\displaylines{\|Q_-(t)-P_-\|_{\frak B}\leqslant {4\cdot\Sp P_-\over
		d_\Pi(B)}\cdot\|C(t)\|_{\frak B},\quad
	\|Q_+(t)-P_+\|_{\frak B}\leqslant {4\cdot\Sp P_+\over d_\Pi(B)}\cdot
		\|C(t)\|_{\frak B},\cr
	\|\dot Q_-\|_{L_1([0,1];\,\frak B)}\leqslant {4\cdot\Sp P_-\over
		d_\Pi^2(B)}\cdot\varkappa(B,C),\quad
	\|\dot Q_+\|_{L_1([0,1];\,\frak B)}\leqslant {4\cdot\Sp P_+\over
		d_\Pi^2(B)}\cdot\varkappa(B,C).
	}
$$
Соответственно, положив
$$
	S(t)\rightleftharpoons Q_0(t)P_0+Q_-(t)P_-+Q_+(t)P_+,
$$
мы получаем оператор-функцию, подчинённую оценкам~\eqref{eq:111},
что ввиду~\eqref{eq:99} гарантирует существование операторов $[S(t)]^{-1}$,
подчинённых оценкам~\eqref{eq:112}. Кроме того, справедливы также
оценки~\eqref{eq:113}. Наконец, будут справедливы равенства
$$
	\displaylines{[S(t)]^{-1}\,[B(t)+C(t)]\,S(t)={}\kern 8truecm\cr
		\kern 2truecm {}=B(t)+[S(t)]^{-1}\cdot
		[Q_0(t)C(t)P_0+Q_-(t)C(t)P_-+Q_+(t)C(t)P_+],}
$$
с учётом оценок
$$
	\leqalignno{\bigl\|[S(t)]^{-1}Q_0(t)-P_0\bigr\|_{\frak B}&=
		\bigl\|[S(t)]^{-1}\cdot [Q_0(t)-P_0]\,(1-P_0)\bigr\|_{\frak B}\cr
		&\leqslant {8\cdot\Sp P_0\over d_\Pi(B)}\cdot\|C(t)\|_{\frak B}
	}
$$
и аналогичных оценок для проекторов $P_-$ и $P_+$ означающие справедливость
оценок~\eqref{eq:114}.
\endproof

Для удобства дальнейшей работы мы сформулируем также следующее тривиальное
утверждение.

\subsubsection\label{prop:3}
{\it Пусть справедлива оценка~\eqref{eq:99}, а оператор-функция
$S\in W_1^1([0,1];\,\frak B)$ удовлетво\-ряет оценкам~\eqref{eq:111}
и~\eqref{eq:112}. Тогда для произвольной оператор-функции $A\in L_1([0,1];\,
\frak B)$ справедлива оценка
$$
	\displaylines{\|S^{-1}AS-\Delta_\Pi\circ A\|_{L_1([0,1];\,\frak B)}
		\leqslant{}\kern 8truecm\cr \kern 2truecm{}\leqslant
		3\cdot\|A-\Delta_\Pi\circ A\|_{L_1([0,1];\,\frak B)}+
		{16\cdot\Sp 1\over d_\Pi(B)}\cdot\int_0^1 \|C(t)\|_{\frak B}
		\cdot\|\Delta_\Pi(A(t))\|_{\frak B}\,dt.
	}
$$
}


\section{Преобразования уравнений}
\subsection
Пусть $\frak A$~--- некоторое направленное множество, а $A_\nu\in L_1([0,1];\,
\frak B)$~--- семейство параметризованных значениями $\nu\in\frak A$
оператор-функций. Основной целью настоящего параграфа является изучение вопроса
о сравнении решения граничной задачи
$$
	\dot x_\nu=A_\nu x_\nu+f,\quad (P_0+P_-)x_\nu(0)+P_+x_\nu(1)=\xi,
	\leqno(\equation)
$$\label{eq:100}%
где $f\in L_1([0,1];\,\frak H)$ и $\xi\in\im P_0$, с решением вспомогательной
граничной задачи
$$
	\dot x_{\nu,0}=A_{\nu,0}x_{\nu,0}+f,\quad
		(P_0+P_-)x_{\nu,0}(0)+P_+x_{\nu,0}(1)=\xi.\leqno(\equation)
$$\label{eq:101}%
Здесь и далее мы используем сокращение $A_{\nu,0}\rightleftharpoons\Delta_\Pi\circ
A_\nu$. Построение преобразований, сводящих задачу~\eqref{eq:100} к подпадающим
под действие утверждения~\ref{prop:0} "<малым"> возмущениям задачи~\eqref{eq:101},
мы будем проводить на основе приближения оператор-функций $A_\nu$ рассмотренными
в предыдущем параграфе \hbox{$\Pi$-кар}\-ка\-са\-ми.

А именно, имеет место следующий факт об асимптотическом поведении решений
параметрических семейств обыкновенных дифференциальных уравнений.

\subsubsection\label{prop:1}
{\it Пусть для некоторой величины $\gamma\in (0,1)$ независимо от выбора значения
$\nu\in\frak A$ и отрезка $[a,b]\subseteq [0,1]$ выполняется оценка
$$
	\int_a^b\omega[(P_0+P_--P_+) A_{\nu,0}(t)]\,dt\leqslant-\ln\gamma.
$$
Пусть также для любого $\varepsilon>0$ найдётся такое значение $\nu_0\in\frak A$,
что для всякого $\nu>\nu_0$ найдётся \hbox{$\Pi$-кар}\-кас $\{B,C\}$ со свойствами
$$
	\varkappa(B,C)<\varepsilon d_\Pi^2(B),\leqno(\equation)
$$\label{eq:151}%
$$
	\abovedisplayskip -1\jot
	\int_0^1\|C(t)\|_{\frak B}^2\,dt<\varepsilon d_\Pi(B),\leqno(\equation)
$$\label{eq:152}%
$$
	\abovedisplayskip -1\jot
	\int_0^1\biggl[\|C(t)\|_{\frak B}\cdot\|A_{\nu,0}(t)-
		B(t)\|_{\frak B}\biggr]\,dt<\varepsilon d_\Pi(B),\leqno(\equation)
$$\label{eq:153}%
$$
	\abovedisplayskip -1\jot
	\|(A_\nu-C)-\Delta_\Pi\circ (A_\nu-C)\|_{L_1([0,1];\,\frak B)}<\varepsilon.
	\leqno(\equation)
$$\label{eq:154}%
Тогда для любого $\varepsilon>0$ найдётся такое $\nu_0\in\frak A$, что для всякого
$\nu>\nu_0$ решения граничных задач~\eqref{eq:100} и~\eqref{eq:101} будут
однозначно определены независимо от выбора $f\in L_1([0,1];\,\frak H)$
и $\xi\in\im P_0$, причём будет также справедлива оценка
$$
	\|x_\nu-x_{\nu,0}\|_{C([0,1];\,\frak H)}\leqslant\varepsilon\cdot
		\bigl[\|x_{\nu,0}\|_{C([0,1];\,\frak H)}+
		\|f\|_{L_1([0,1];\,\frak H)}\bigr].\leqno(\equation)
$$\label{eq:301}%
}%
\proof
Подстановкой $x_\nu=Sy$ граничная задача~\eqref{eq:100} приводится
к эквивалентному виду
$$
	\dot y=\bigl\{A_{\nu,0}+[S^{-1}A_\nu S-A_{\nu,0}-S^{-1}\dot S]\bigr\}\,y
		+S^{-1}f,\quad (P_0+P_-)y(0)+P_+y(1)=\xi.
$$
При этом, с учётом сделанных предположений и утверждений~\ref{prop:kar}
и~\ref{prop:3}, имеют место оценки
$$
	\|S^{-1}A_\nu S-A_{\nu,0}-S^{-1}\dot S\|_{L_1([0,1];\,
		\frak B)}\leqslant (48\cdot\Sp 1+3)\,\varepsilon.
$$
Соответственно, в случае, когда значение $\theta\rightleftharpoons
(48\cdot\Sp 1+3)\,\varepsilon\gamma^{-1}$ лежит на интервале $(0,1)$, искомое
решение $y\in W_1^1([0,1];\,\frak H)$ однозначно определено и подчиняется оценке
$$
	\|y-z\|_{C([0,1];\,\frak H)}\leqslant
		{\theta\over 1-\theta}\,\|z\|_{C([0,1];\,\frak H)},
$$
где $z\in W_1^1([0,1];\,\frak H)$ есть решение граничной задачи
$$
	\dot z=A_{\nu,0}z+S^{-1}f,\quad (P_0+P_-)z(0)+P_+z(1)=\xi.
$$
Учёт оценок
$$
	\leqalignno{\|z-x_{\nu,0}\|_{C([0,1];\,\frak H)}&\leqslant\gamma^{-1}\,
		\|S^{-1}-1\|_{C([0,1];\,\frak B)}\,\|f\|_{L_1([0,1];\,\frak H)},\cr
	\|S(t)-1\|_{\frak B}&\leqslant {4\cdot\Sp 1\over d_\Pi(B)}
		\cdot\|C(t)\|_{\frak B}\cr
		&\leqslant {2\cdot\Sp 1\over d_\Pi(B)}\cdot
		\int_0^1 \|\dot C(s)\|_{\frak B}\,ds\cr
		&\leqslant {2\cdot\Sp 1\over d_\Pi^2(B)}\cdot\varkappa(B,C)\cr
		&<2\varepsilon\cdot\Sp 1
	}
$$
позволяет теперь завершить доказательство.
\endproof

Утверждение~\ref{prop:1} очевидным образом имеет следующий двойственный
относительно отражения отрезка $[0,1]$ вариант.

\subsubsection\label{prop:1a}
{\it Пусть для некоторой величины $\gamma\in (0,1)$ независимо от выбора значения
$\nu\in\frak A$ и отрезка $[a,b]\subseteq [0,1]$ выполняется оценка
$$
	\int_a^b\omega[(P_--P_0-P_+) A_{\nu,0}(t)]\,dt\leqslant-\ln\gamma.
$$
Пусть также для любого $\varepsilon>0$ найдётся такое значение $\nu_0\in\frak A$,
что для всякого $\nu>\nu_0$ найдётся \hbox{$\Pi$-кар}\-кас $\{B,C\}$
со свойствами~\eqref{eq:151}--\eqref{eq:154}. Тогда для любого $\varepsilon>0$
найдётся такое $\nu_0\in\frak A$, что для всякого $\nu>\nu_0$ решения граничных
задач
$$
	\displaylines{\dot x_\nu=A_\nu x_\nu+f,\quad P_-x_\nu(0)+
		(P_0+P_+)x_\nu(1)=\xi,\cr
	\dot x_{\nu,0}=A_{\nu,0}\,x_{\nu,0}+f,\quad
		P_-x_{\nu,0}(0)+(P_0+P_+)x_{\nu,0}(1)=\xi
	}
$$
будут однозначно определены независимо от выбора $f\in L_1([0,1];\,\frak H)$
и $\xi\in\im P_0$, причём будет также справедлива оценка~\eqref{eq:301}.
}%

\subsection
Следует заметить, что после того, как справедливость оценок вида~\eqref{eq:301}
установлена, эти оценки допускают несложное дальнейшее уточнение. В частности,
справедливо следующее утверждение, аналоги которого широко используются
в работах~\hbox{[\cite{SSh:2020}, \cite{Sh:2021}]}.

\subsubsection
{\it Пусть выполнены предположения утверждения~\ref{prop:1}. Тогда для любого
$\varepsilon>0$ найдётся такое $\nu_0\in\frak A$, что для всякого $\nu>\nu_0$
решения граничных задач~\eqref{eq:100} и~\eqref{eq:101} будут независимо от выбора
$f\in L_1([0,1];\,\frak H)$ и $\xi\in\im P_0$ подчинены оценке
$$
	\|x_\nu-x_{\nu,0}\|_{C([0,1];\,\frak H)}\leqslant
		\left[1+(\gamma^{-1}+\varepsilon)\cdot
		\|A_\nu-A_{\nu,0}\|_{L_1([0,1];\,\frak H)}\right]\cdot
		\|y_\nu\|_{C([0,1];\,\frak H)},
$$
где функция $y_\nu\in W_1^1([0,1];\,\frak H)$ представляет собой решение
граничной задачи
$$
	\dot y_\nu=A_{\nu,0}y_\nu+[A_\nu-A_{\nu,0}]\,x_{\nu,0},\qquad
		(P_0+P_-)y_\nu(0)+P_+y_\nu(1)=0.
$$
}%
\proof
Решение $z\in W_1^1([0,1];\,\frak H)$ граничной задачи
$$
	\dot z=A_\nu z+[A_\nu-A_{\nu,0}]\,y_\nu,\qquad (P_0+P_-)z(0)+P_+z(1)=0
$$
удовлетворяет равенству $x_\nu-x_{\nu,0}=y_\nu+z$. При этом, согласно
утверждению~\ref{prop:1}, при достаточно больших $\nu\in\frak A$ имеет место оценка
$$
	\|z\|_{C([0,1];\,\frak H)}\leqslant\left[(1+\varepsilon)\,\gamma^{-1}+
		\varepsilon\right]\cdot
		\|A_\nu-A_{\nu,0}\|_{L_1([0,1];\,\frak H)}\cdot
		\|y_\nu\|_{C([0,1];\,\frak H)}.
$$
Это и означает справедливость доказываемого утверждения.
\endproof

\subsection
Рассмотрим теперь ситуацию, когда задана полная система $\{P_k\}_{k=1}^m$ ненулевых
ортопроекторов на попарно ортогональные подпространства, а оператор-функции
$A_\nu$ допускают представление
$$
	A_\nu=\sum_{k=1}^m d_\nu h_{\nu,k}\cdot P_k+V_\nu.
$$
Здесь $\{d_\nu\}_{\nu\in\frak A}$ есть некоторая бесконечно большая числовая
направленность, а функции $h_{\nu,k}\in L_1[0,1]$ и $V_\nu\in L_1([0,1];\,\frak B)$
предполагаются подчинёнными следующим условиям:

\subsubsection\label{sup:11}
Для любого $\varepsilon>0$ найдутся такое значение $\nu_0\in\frak A$ и шар $I$
пространства $W_1^1[0,1]$, что для всяких $\nu>\nu_0$ и $k\in 1\,.\,.\,m$
расстояние между функцией $h_{\nu,k}$ и образом вложения шара $I$ в пространство
$L_1[0,1]$ не будет превосходить выбранного значения $\varepsilon>0$.

\subsubsection\label{sup:12}
Для любого $\varepsilon>0$ найдутся такое значение $\nu_0\in\frak A$ и шар $I$
пространства $\Wo_1^1([0,1];\,\frak B)$, что для всякого $\nu>\nu_0$ расстояние
между оператор-функцией $V_\nu$ и образом вложения шара $I$ в пространство
$L_1([0,1];\,\frak B)$ не будет превосходить выбранного значения $\varepsilon>0$.

\subsubsection\label{sup:13}
Для любого $\nu\in\frak A$ при почти всех $t\in [0,1]$ независимо от выбора
значений $k,l\in 1\,.\,.\,m$ со свойством $k<l$ выполняются оценки
$$
	|h_{\nu,k}(t)-h_{\nu,l}(t)|\geqslant 1.
$$

\subsubsection\label{sup:14}
Существует величина $\gamma\in (0,1)$, для которой независимо от выбора отрезка
$[a,b]\subseteq [0,1]$, значения $\nu\in\frak A$ и значений $k,l\in 1\,.\,.\,m$
со свойством $k<l$ выполняются оценки
$$
	\int_a^b \Re[d_\nu h_{\nu,k}(t)-d_\nu h_{\nu,l}(t)]\,dt\leqslant-\ln\gamma.
$$

\medskip
Следует заметить, что типичной в приложениях (см., например,
\hbox{[\cite{SSh:2020}]}) является ситуация, когда функции $h_{\nu,k}$ и $V_\nu$
не зависят от выбора значения $\nu\in\frak A$. При этом условия~\ref{sup:11}
и~\ref{sup:12} выполняются автоматически ввиду плотности вложений пространств
$W_1^1[0,1]$ и $\Wo_1^1([0,1];\,\frak B)$ в пространства $L_1[0,1]$
и $L_1([0,1];\,\frak B)$, соответственно.

Из утверждения~\ref{prop:1} вытекает справедливость следующего утверждения.

\subsubsection\label{prop:2}
{\it Пусть выполнены условия~\ref{sup:11}--\ref{sup:14}. Пусть также зафиксированы
значение $k\in 1\,.\,.\,m$ и ненулевой вектор $\xi\in\im P_k$. Пусть, наконец,
положено
$$
	P_0\rightleftharpoons P_k,\qquad P_-\rightleftharpoons\sum_{l=1}^{k-1} P_l,
	\qquad P_+\rightleftharpoons\sum_{l=k+1}^m P_l.\leqno(\equation)
$$\label{eq:1}%
Тогда решения граничных задач
$$
	\displaylines{\dot x_\nu=A_\nu x_\nu,\quad (P_0+P_-)x_\nu(0)+
		P_+x_\nu(1)=\xi,\cr
	\dot x_{\nu,0}=A_{\nu,0}x_{\nu,0},\quad
		(P_0+P_-)x_{\nu,0}(0)+P_+x_{\nu,0}(1)=\xi
	}
$$
однозначно определены при всех достаточно больших $\nu\in\frak A$ и подчинены
асимптотической оценке
$$
	\sup_{t\in [0,1]}{\|x_\nu(t)-x_{\nu,0}(t)\|_{\frak H}\over
		\|x_{\nu,0}(t)\|_{\frak H}}\to 0.
$$
}%
\proof
Для доказательства искомого достаточно установить приложи\-мость заключения
утверждения~\ref{prop:1} к вспомогательной системе
$$
	\dot x_\nu=[A_\nu-d_\nu h_{\nu,k}]x_\nu.\leqno(\equation)
$$\label{eq:302}%

Зафиксируем произвольное значение $\varepsilon>0$, а также такие $\nu_0\in\frak A$
и $\alpha,\beta>0$, что при всяком $\nu>\nu_0$ найдутся оператор-функция
$C\in\Wo_1^1([0,1];\,\frak B)$ со свойствами
$$
	\|V_\nu-C\|_{L_1([0,1];\,\frak B)}<\varepsilon/2,\quad
	\|\dot C\|_{L_1([0,1];\,\frak B)}<\alpha,
$$
а также набор $\{g_l\}_{l=1}^m$ функций класса $W_1^1[0,1]$ со свойствами
$$
	\displaylines{\|g_l-h_{\nu,l}\|_{L_1[0,1]}<{\varepsilon\over 4\alpha m},
	\qquad \|\dot g_l\|_{L_1[0,1]}<\beta,\cr
	|g_l(t)-g_s(t)|\geqslant 1/2\quad\hbox{при $l\neq s$.}
	}
$$
Тогда для оператор-функции вида
$$
	B\rightleftharpoons\sum_{l=1}^m d_\nu\cdot[g_l-g_k]\cdot P_l
$$
будут справедливы оценки
$$
	\displaylines{d_\Pi(B)\geqslant |d_\nu|/2,\qquad
	\|\dot B\|_{L_1([0,1];\,\frak B)}\leqslant 2\beta m\,|d_\nu|,\cr
	\|\Delta_\Pi\circ [A_\nu-d_\nu h_{\nu,k}]-B\|_{L_1([0,1];\,\frak B)}<
		{\varepsilon |d_\nu|\over 2\alpha}+
		\|V_\nu\|_{L_1([0,1];\,\frak B)}.
	}
$$
Оценка~\eqref{eq:154} для системы~\eqref{eq:302} при этом выполняется
автоматически. Кроме того, ввиду вытекающей из сделанных предположений оценки
$\|C\|_{C([0,1];\,\frak B)}<\alpha/2$, при больших значениях параметра $d_\nu$
пара $\{B,C\}$ будет являться \hbox{$\Pi$-кар}\-ка\-сом, удовлетворяющим
оценкам~\eqref{eq:151}, \eqref{eq:152} и~\eqref{eq:153}. Тем самым, доказываемое
утверждение справедливо.
\endproof

Аналогичным образом, из утверждения~\ref{prop:1a} вытекает справедливость
следующего утверждения.

\subsubsection\label{prop:2a}
{\it Пусть выполнены условия~\ref{sup:11}--\ref{sup:14}. Пусть также зафиксированы
значение $k\in 1\,.\,.\,m$ и ненулевой вектор $\xi\in\im P_k$. Пусть, наконец,
положено~\eqref{eq:1}. Тогда решения граничных задач
$$
	\displaylines{\dot x_\nu=A_\nu x_\nu,\quad
		P_-x_\nu(0)+(P_0+P_+)x_\nu(1)=\xi,\cr
	\dot x_{\nu,0}=A_{\nu,0}\,x_{\nu,0},\quad
		P_-x_{\nu,0}(0)+(P_0+P_+)x_{\nu,0}(1)=\xi
	}
$$
однозначно определены при всех достаточно больших $\nu\in\frak A$ и подчинены
асимптотической оценке
$$
	\sup_{t\in [0,1]}{\|x_\nu(t)-x_{\nu,0}(t)\|_{\frak H}\over
		\|x_{\nu,0}(t)\|_{\frak H}}\to 0.
$$
}%

\subsection
Чтобы непосредственно подпасть под действие утверждения~\ref{prop:1} и его
вариантов, система уравнений должна быть в некотором смысле "<почти диагональной">.
В типичном случае произвольная система вида $\dot x=(A+V)x$, где оператор-функция
$A$ абсолютно непрерывна, может быть приведена к такой "<почти диагональной"> форме
заменой $x=Sy$, где $S(t)$~--- преобразование, диагонализирующее оператор $A(t)$.
Конкретный пример задачи, легко разбираемой на основе такой схемы, будет указан
далее.


\section{Пример}
\subsection
Пусть $\frak H=\Bbb C^n$, где $n\geqslant 2$, и пусть изучаемая система
имеет вид
$$
	\dot x_\lambda=D_\lambda x_\lambda,\qquad
	(D_\lambda)_{kl}=\cases{p_k&при $l=k+1$,\cr
		(\lambda+\zeta)^n p_n-q_{n1}&при $k=n$, $l=1$,\cr
		-q_{kl}&при $k+1-n<l\leqslant k$,\cr 0&иначе.}
$$
Здесь $\lambda\in\Bbb C\setminus\{0\}$~--- параметр, $\zeta\in\Bbb C$~--- некоторое
фиксированное число, функции $p_k\in W_1^1[0,1]$ принимают положительные значения,
а функции $q_{kl}\in L_1[0,1]$ выбраны произвольным образом. Такого рода системы,
как известно, возникают при изучении обыкновенных дифференциальных операторов
различных порядков \hbox{[\cite{Na:1969}: \S$\,$16]}, в том числе
с коэффициентами-обобщёнными функциями \hbox{[\cite{SSh:2003}, \cite{Vl:2004},
\cite{Vl:2017}, \cite{SSh:2020}]}.

Оператор-функция $D_\lambda$ обладает "<ведущей частью">
$$
	(D_{\lambda,0})_{kl}=\cases{p_k&при $l=k+1$,\cr
		\lambda^n p_n&при $k=n$, $l=1$,\cr 0&иначе,}
$$
имеющей диагонализирующее разложение
$$
	[S_\lambda(t)]^{-1}D_{\lambda,0}(t)S_\lambda(t)=A_\lambda(t),\qquad
		[A_\lambda(t)]_{kk}\rightleftharpoons
		\left[\lambda\cdot\left(\prod_{m=1}^n p_m(t)\right)^{1/n}\right]
		\cdot e^{2\pi ki/n},
$$
где положено
$$
	[S_\lambda(t)]_{lk}=[A_\lambda(t)]_{kk}^{l-1}\cdot
		\prod_{m=1}^{l-1} p^{-1}_m(t),\qquad
	[S_\lambda(t)]^{-1}_{kl}={[A_\lambda(t)]_{kk}^{1-l}\over n}
		\cdot\prod_{m=1}^{l-1}p_m(t).
$$

Оператор-функции $S_\lambda^{-1}\dot S_\lambda$ не зависят от выбора значения
$\lambda\in\Bbb C\setminus\{0\}$, причём для всех $k\in 1\,.\,.\,n$ выполняются
равенства
$$
	\bigl[[S_\lambda(t)]^{-1}\dot S_\lambda(t)\bigr]_{kk}=
		\sum_{m=1}^n{2m-n-1\over 2n}\cdot {d\ln p_m\over dt}.
$$
Далее, для операторов $G,F_m\in\frak B$, где $m\in 1\,.\,.\,n$, с матричными
элементами
$$
	G_{kl}\rightleftharpoons e^{2\pi ki/n},\qquad
	(F_m)_{kl}\rightleftharpoons n^{-1}\cdot e^{2\pi (l-k)(m-1)i/n}
$$
имеет место асимптотическое соотношение
$$
	\biggl\|S_\lambda^{-1}\,[D_\lambda-D_{\lambda,0}] S_\lambda-
		\zeta\cdot\left(\prod_{m=1}^n p_m\right)^{1/n}\cdot G+
		\sum_{m=1}^n q_{mm}\cdot F_m\biggr\|_{L_1([0,1];\,\frak B)}
		=O(|\lambda|^{-1}).
$$
Наконец, при любом выборе значения $m\in 1\,.\,.\,2n$ для сектора
$$
	\Omega_m\rightleftharpoons\left\{\lambda\in\Bbb C\setminus\{0\}\;:\;
		\Arg\lambda\in\left({\pi\cdot(m-1)\over n},
		{\pi m\over n}\right)\right\}
$$
может быть указана такая перестановка $\tau_m$ списка $1\,.\,.\,n$,
что соответствующий полный набор $\{P_k\}_{k=1}^n$ диагональных ортопроекторов вида
$$
	(P_k)_{ll}=\cases{1&при $l=\tau_m(k)$,\cr 0&иначе}
$$
вместе с системой
$$
	\dot y_\lambda=[S_\lambda^{-1}D_\lambda S_\lambda-
		S_\lambda^{-1}\dot S_\lambda]\,y_\lambda\leqno(\equation)
$$\label{eq:303}%
будут удовлетворять условиям утверждений~\ref{prop:2} и~\ref{prop:2a}.
Соответственно, при $|\lambda|\gg 0$ у системы~\eqref{eq:303} существуют
аналитически зависящие от параметра $\lambda\in\Omega_m$ решения вида
$$
	\displaylines{y_{m,k,\lambda}(t)=\left(\prod_{l=1}^n\,
		[p_l(t)]^{n-2l+1\over 2n}\right)\times{}\kern 9truecm\cr
		\kern 2truecm{}\times\exp\left(\int_0^t\left[(\lambda+\zeta)\cdot
		\left(\prod_{l=1}^n p_n(s)\right)^{1/n}\cdot e^{2\pi ki/n}-
		{1\over n}\,\sum_{l=1}^n q_{ll}(s)\right]\,ds\right)
		\cdot[\xi_k+o(1)],
	}
$$
где векторы $\xi_k\in\frak H$ имеют вид $(\xi_k)_l=\delta_{kl}$.

\vskip 0.4cm
\eightpoint\rm
{\leftskip 0cm\rightskip 0cm plus 1fill\parindent 0cm
\bf Литература\par\penalty 20000}\vskip 0.4cm\penalty 20000
\bibitem{SSh:2020} Савчук~А.$\,$М., Шкаликов~А.$\,$А. {\it Асимптотический анализ
решений обыкновенных дифференциальных уравнений
с коэффициентами-распределениями}~// Матем. сборник.~--- 2020.~--- Т.~211, \No~11.~---
С.~129--166.
\bibitem{Sh:2021} Шкаликов~А.$\,$А. {\it Регулярные спектральные задачи для систем
обыкновенных дифференциальных уравнений первого порядка}~// Успехи матем. наук.~---
2021.~--- Т.~76, \No~5.~--- С.~203--204.
\bibitem{Na:1969} Наймарк~М.$\,$А. {\it Линейные дифференциальные операторы\/}.
М.:~Наука, 1969.
\bibitem{P:1974} Понтрягин~Л.$\,$С. {\it Обыкновенные дифференциальные уравнения\/},
изд.~4. М.:~Наука, 1974.
\bibitem{SSh:2003} Савчук~А.$\,$М., Шкаликов~А.$\,$А.
{\it Операторы Штурма--Лиувилля с потенциалами-распределениями}~// Труды~ММО.~---
2003.~--- Т.~64.~--- С.~159--212.
\bibitem{Vl:2004} Владимиров~А.$\,$А. {\it О сходимости последовательностей
обыкновенных дифференциальных операторов}~// Матем.~заметки.~--- 2004.~---
Т.~75, \No~6.~--- С.~941--943.
\bibitem{Vl:2017} Владимиров~А.$\,$А. {\it Об одном подходе к определению
сингулярных дифференциальных операторов}~// arXiv:1701.08017.~--- 2017.~--- 11~стр.
\bye